\documentclass[a4paper,12pt,leqno]{article}
\usepackage{latexsym}
\usepackage[all]{xy}

\usepackage{amssymb} 
\usepackage{amsmath} 
\usepackage{theorem}

\def\Z{{\mathbb{Z}}}
\def\R{{\mathbb{R}}}
\def\K{{\mathbb{K}}}
\def\A{{\mathcal{A}}}
\def\B{{\mathcal{B}}}

\DeclareMathOperator{\coker}{coker}

\DeclareMathOperator{\Der}{Der}

\numberwithin{equation}{section}

\newcommand{\owari}{\hfill$\square$}

\theoremstyle{break}
\newtheorem{theorem}{Theorem}[section]
\newtheorem{prop}[theorem]{Proposition}
\newtheorem{cor}[theorem]{Corollary}
\newtheorem{lemma}[theorem]{Lemma}
\newtheorem{define}[theorem]{Definition}
\newtheorem{rem}[theorem]{Remark}

\newtheorem{example}[theorem]{Example}
\newtheorem{conj}[theorem]{Conjecture}

\newcommand{\xgraphAvertex}[1][****]{
\xgraphAVertex #1
}
\newcommand{\xgraphAVertex}[4]{
\if#3o\put(0,0){\circle{4}}\fi%
\if#3*\put(0,0){\circle*{4}}\fi%
\if#3.\put(0,0){\circle*{2.4}}\fi%
\if#4o\put(30,0){\circle{4}}\fi%
\if#4*\put(30,0){\circle*{4}}\fi%
\if#4.\put(30,0){\circle*{2.4}}\fi%
\if#2o\put(0,30){\circle{4}}\fi%
\if#2*\put(0,30){\circle*{4}}\fi%
\if#2.\put(0,30){\circle*{2.4}}\fi%
\if#1o\put(30,30){\circle{4}}\fi%
\if#1*\put(30,30){\circle*{4}}\fi%
\if#1.\put(30,30){\circle*{2.4}}\fi%
}
\newcommand{\xgraphA}[6]{
\if#1+\put(30,30){\line(-1,0){30}}\fi
\if#1.\qbezier[7](30,30)(15,30)(0,30)\fi
\if#1-\put(30,31){\line(-1,0){30}}\put(30,29){\line(-1,0){30}}\fi
\if#2+\put(0,0){\line(1,1){30}}\fi 
\if#2.\qbezier[10](0,0)(15,15)(30,30)\fi 
\if#2-\put(-0.7,0.7){\line(1,1){30}}\put(0.7,-0.7){\line(1,1){30}}\fi 
\if#3+\put(30,30){\line(0,-1){30}}\fi
\if#3.\qbezier[7](30,0)(30,15)(30,30)\fi 
\if#3-\put(29,30){\line(0,-1){30}}\put(31,30){\line(0,-1){30}}\fi
\if#4+\put(0,0){\line(0,1){30}}\fi 
\if#4.\qbezier[7](0,0)(0,15)(0,30)\fi 
\if#4-\put(-1,0){\line(0,1){30}}\put(1,0){\line(0,1){30}}\fi 
\if#5+\put(0,30){\line(1,-1){30}}\fi 
\if#5.\qbezier[10](30,0)(15,15)(0,30)\fi 
\if#5-\put(-0.7,29.3){\line(1,-1){30}}\put(0.7,30.7){\line(1,-1){30}}\fi 
\if#6+\put(0,0){\line(1,0){30}}\fi 
\if#6.\qbezier[7](0,0)(15,0)(30,0)\fi 
\if#6-\put(0,1){\line(1,0){30}}\put(0,-1){\line(1,0){30}}\fi 
\xgraphAvertex}

\title{Exponents of $2$-multiarrangements and freeness of $3$-arrangements
}
\author{Takuro Abe\thanks{
Department of Mathematics, Kyoto University, 
Kitashirakawa-Oiwake-Cho, Sakyo-Ku, 
Kyoto 606--8502, Japan.
email:abetaku@math.kyoto-u.ac.jp.}}
\date{\today}

\begin{document}

\maketitle

\begin{abstract}
We give the upper bound of differences 
of exponents for balanced $2$-multiarrangements in terms of 
the cardinality of hyperplanes. 
Also, we give a shift isomorphism of $2$-multiarrangements 
like Coxeter arrangements when the difference of exponents is maximum. 
As an application, a sufficient numerical and combinatorial condition for 
$3$-arrangements to be free is given. 
\end{abstract}
\setcounter{section}{-1}

\section{Introduction}
Let $V$ be an $\ell$-dimensional vector space over a field $\K$ of characteristic zero, 
$S=\K[x_1,\ldots,x_\ell]$ the coordinate ring and 
$\Der(S)=\oplus_{i=1}^\ell S \cdot \partial_{x_i}$ the module of 
$S$-regular derivations. A \textbf{hyperplane arrangement} $\A$ is 
a finite collection of hyperplanes in $V$. In this article $\A$ is assumed to 
consist of linear hyperplanes unless otherwise specified. Such an arrangement 
is called \textbf{central}. 
For each $H \in \A$ let us fix a linear form $\alpha_H \in V^*$ such that 
$\ker \alpha_H=H$. A function 
$m:\A \rightarrow \Z_{>0}$ is called a \textbf{multiplicity} and 
a pair $(\A,m)$ is a \textbf{multiarrangement}. 
Then we can define the \textbf{logarithmic derivation module} 
$D(\A,m)$ by 
$$
D(\A,m):=\{\theta \in \Der(S) \mid 
\theta(\alpha_H) \in S \cdot \alpha_H^{m(H)}\ (\forall H \in \A)\}.
$$
$D(\A,m)$ is a reflexive module in general. 
When $D(\A,m)$ is a free module of rank $\ell$, we say that 
$(\A,m)$ is \textbf{free} and for the homogeneous basis $\theta_1,\ldots,
\theta_\ell$, we define 
$$
\exp(\A,m):=(\deg \theta_1,\ldots,\deg \theta_\ell),
$$
where $\deg \theta:=\deg \theta(\alpha)$ for some $\alpha \in V^*$ such that 
$\theta(\alpha)\neq 0$. 
When $m \equiv 1$ a multiarrangement $(\A,1)$ is the same as an arrangement, which 
is sometimes called a \textbf{simple arrangement} and 
$D(\A,1)=:D(\A)$. An \textbf{$\ell$-arrangement} is 
that in $V \simeq \K^\ell$. 

The freeness of an arrangement $\A$ has been studied by a lot of mathematicians for a long 
time. Actually it is very difficult to determine whether a given arrangement is free or not. 
For example, whether the freeness of simple arrangements depends only on 
the combinatorics of arrangements or not has been unsolved for a long time, which is called the 
\textbf{Terao conjecture} and still open. 
Recently, new freeness criterions were found by Yoshinaga in \cite{Y2} and \cite{Y3} in terms of 
restricted multiarrangements. 
Hence it has become important to study 
the freeness of multiarrangements. In particular, by \cite{Y3}, to solve the Terao conjecture of 
$3$-arrangements, to determine exponents of $2$-multiarrangements is essential. 

Since $2$-multiarrangements are free, we can always define the exponents. 
However, contrary to the simple arrangement case, the behavior of $\exp(\A,m)$ is 
complicated. One of the approaches to understand it 
is \cite{AN} in which a multiplicity lattice is introduced and studied. 
The aim of this article is the further analysis of the theory of 
multiplicity lattices, introduction of a generalized Euler derivations called 
$(\A,m)$-Euler derivations, and 
apply them to the freeness problem of $3$-arrangements as desired. Let us 
explain these in details below.

For a $2$-multiarrangement $(\A,m)$ with $\exp(\A,m)=(d_1,d_2)$, let us define 
$$
\Delta(m):=|d_1-d_2|.
$$
%
We say that a multiplicity $m$ is \textbf{balanced} if 
$$
m(K) \le \sum_{H \in \A \setminus \{K\}} m(H)\ (\forall K \in \A).
$$
In \cite{AN} the structure of multiplicity lattices was considered and 
studied by using $\Delta$. The theory constructed there will be 
used to prove results in this article. See \cite{AN}, Lemma \ref{one}, 
Theorem \ref{str} and section one for details. 

Note that, if $m$ is not balanced, then $\exp(\A,m)$ can be easily computed,  
see Proposition \ref{nb}. Hence for the Terao conjecture of $3$-arrangements, 
we have to study the exponents of balanced $2$-multiarrangements. Since 
$d_1+d_2=|m|$ when $\exp(\A,m)=(d_1,d_2)$, to know exponents is 
equivalent to know $\Delta(m)$. 
Then it is a natural question to ask 
for a $2$-multiarrangement $(\A,m)$, is there any upper bound 
of $\Delta(m)$ when $m$ is balanced? Experimental computations imply that 
$|\A|-2$ might be the upper bound. In fact, $\Delta(m)=h-2$ when $m \equiv 1$. 
The first main result in this article is to prove that it is in fact the 
strict upper bound.

\begin{theorem}
Let $\A$ be a $2$-arrangement with $|\A|=h>2$. If $m:\A \rightarrow \Z_{>0}$ is 
balanced, then $\Delta(m) \le h-2$. 
\label{limit}
\end{theorem}

For the proof, we use results in \cite{AN} and the affine connection $\nabla$. 
Then it is an interesting problem to ask whether 
there are some special properties if $m$ is balanced and $\Delta(m)=|\A|-2$. 
When $m \equiv 1$, this condition is satisfied. Let us agree that 
the \textbf{lower degree basis} $\theta$ for $D(\A,m)$ 
is the homogeneous derivation $\theta$ such that 
$\{\theta,\varphi\}$ is an $S$-basis for $D(\A,m)$ and that 
$\deg \theta \le \deg \varphi$. Then the lower degree basis 
for $D(\A)$ is the Euler derivation, which is apparently special.  
The 
answer is 
interesting. 

\begin{theorem}
Let $\A$ be an arrangement in $\K^2$ with $|\A|=h>2$ and 
$m_0:\A \rightarrow \Z_{> 0}$ be 
a balanced multiplicity such that 
$\Delta(m_0)=h-2$.
Assume that one of the following two holds:
\begin{itemize}
\item[(1)]
$|\A|=h=3$ and $m_0-1$ is balanced, or 
\item[(2)]
$|\A|=h \ge 4$.
\end{itemize}
Then the lower degree basis $\theta_0$ for $D(\A,m_0)$ gives rise to an isomorphism 
$$
\Phi_0:D(\A,m) \rightarrow D(\A,m_0+m-1)
$$
defined by 
$$
\Phi_0(\theta):=\nabla_\theta \theta_0
$$
with 
$$
m:\A \rightarrow \{+1,0\}.
$$
\label{univ}
\end{theorem}

The isomorphism $\Phi_0$ is first 
introduced in \cite{Y0}, generalized in \cite{AY2} and 
\cite{A4} all for Coxeter multiarrangements. In these papers, 
the invariant theory of Coxeter groups and the existence of 
the primitive derivation played important roles. On the contrary, Theorem \ref{univ} 
do not need them and the same statement can be true for all $2$-arrangements.  

As an application of Theorem \ref{limit} a 
freeness condition for $3$-arrangement is also given. It is well-known that 
when $\A$ is free with $\exp(\A)=(d_1,\ldots,d_\ell)$ the characteristic polynomial 
$\chi(\A,t)$ 
splits into 
\begin{equation}
\chi(\A,t)=\prod_{i=1}^\ell (t-d_i).
\label{tf}
\end{equation}
The formula (\ref{tf}) is the famous Terao's factorization theorem proved in \cite{Tf}. 
However, the converse does not hold. Theorem \ref{limit} combined with 
the result in \cite{Y3} gives the converse for some cases. In other words, 
if a $3$-arrangement has a splitting characteristic polynomial with 
certain exponents, then it is free as follows. 

\begin{theorem}
Let $\overline{\A}$ be an affine $2$-arrangement, $c\overline{\A}$ its coning 
with the 
infinite hyperplane $H_0 \in c \overline{\A}$. Put $|\overline{\A}|=k$ and 
$\chi(\overline{\A},t)=(t^2-kt+c_2)$. Also, let $(\A,m)$ be 
the Ziegler restriction (see Definition \ref{zr}) of $c \overline{\A}$ onto $H_0$ with $|\A|=h>2$. 

If $(\A,m)$ is balanced and 
$\chi(\overline{\A},t)=(t-d)(t-d-h+2)$ or $\chi(\overline{\A},t)=(t-d)(t-d-h+3)$ for 
some integer $d$, then 
$c\overline{\A}$ is free.
\label{fc}
\end{theorem}

Note that every central $3$-arrangement can be obtained as 
the coning of a certain affine $2$-arrangement. Hence 
Theorem \ref{fc} says that only the combinatorics determines the freeness of 
some $3$-arrangements. 
More explicitly, if we define the category of $3$-arrangements 
$PB_3$ 
which consists of $\A$ 
such that every Ziegler restriction is balanced, and there exists $H_0 \in \A$ such that 
$\chi(\A,t)=(t-1)(t-d)(t-d')$ for some $d,d' \in \Z$ with 
$|d-d'| \ge h-3$, where 
$h=|\A \cap H_0|$. Then we have the following:

\begin{cor}
The Terao conjecture \ref{TC} is true in $PB_3$.
\label{TC3}
\end{cor}

Since the Terao conjecture is true for non-balanced $3$-arrangements (see Proposition 
\ref{nb} or \cite{WY}), Corollary \ref{TC3} is the first step for the 
Terao conjecture of balanced $3$-arrangements. 

The organization of this article is as follows. In section one we introduce some 
notions and results which will be used in this article. In section two we prove 
Theorem \ref{limit}. In section three we prove Theorem \ref{univ}. In section four 
we prove Theorem \ref{fc} and Corollary \ref{TC3}. 
Also we give several applications of Theorem \ref{limit} and 
examples of free $3$-arrangements.  
\medskip

\noindent
\textbf{Acknowledgements}. The author 
thanks Yasuhide Numata for his reading the first draft of this article and pointing 
out a mistake. 

\section{Preliminaries}
In this section let us summarize results and definitions which will be 
used in this article. For a general reference, see \cite{OT}. 
We use the notation in the introduction. For 
an affine $\ell$-arrangement $\A$, the \textbf{coning} $c\A$ of $\A$ is 
an $(\ell+1)$-arrangement obtained by replacing 
$
\{\alpha=k\} \in \A\ (\alpha\in V^*,\ k \in \K)$
by 
$
\{\alpha=kz_{\infty}\} \in c\A
$
combined with the \textbf{infinite hyperplane} $H_0=\{z_{\infty}=0\} \in c\A$.
For a simple arrangement $\A$ define the \textbf{intersection lattice} $L(\A)$ by 
$$
L(\A):=\{ \cap_{H \in \B} H \mid \B \subset \A\}.
$$
This is a partially ordered set with the reverse inclusion order and 
the unique minimum element $V \in L(\A)$. The \textbf{M\"obius function} $\mu$ 
on $L(\A)$ is defined by $\mu(V)=1$ and by 
$$
\mu(X):=-\sum_{X \subsetneq Y} \mu(Y)\ (X \in L(\A) \setminus \{V\}).
$$
A \textbf{characteristic polynomial} $\chi(\A,t)$ is defined by 
$$
\chi(\A,t):=\sum_{X \in L(\A)} \mu(X) t^{\dim X}.
$$
The following is one of the most important problems among 
the arrangement theory. 

\begin{conj}[Terao]
The freeness of a simple arrangement $\A$ depends only on 
its intersection lattice $L(\A)$.
\label{TC}
\end{conj}

For a multiarrangement $(\A,m)$, put 
$$
|m|:=\sum_{H \in \A} m(H).
$$
The following is the most fundamental 
result in the free arrangement theory. For the proof, see \cite{OT} and \cite{Z}. 

\begin{theorem}[Saito's criterion]
Let $\theta_1,\ldots,\theta_\ell$ be homogeneous derivations in $D(\A,m)$. 
Then $\A$ is free with basis $\{\theta_1,\ldots,\theta_\ell\}$ if and only if $\{
\theta_1,\ldots,\theta_\ell\}$ is $S$-independent and $\sum_{i=1}^\ell \deg \theta_i=|m|$.
\label{Saito}
\end{theorem}

To use Yoshinaga's freeness criterion, we often use the Ziegler restriction. 

\begin{define}[Ziegler restriction]
Let $\A$ be a simple arrangement and fix $H_0 \in \A$. A \textbf{Ziegler restriction} 
$(\A'',m_0)$ of $\A$ with respect to $H_0$ is defined by 
\begin{eqnarray*}
\A'':&=&\{H \cap H_0 \mid H \in \A \setminus \{H_0\}\},\\
m_0(K):&=&| \{H \in \A \setminus \{H_0\} \mid H \cap H_0 =K\}|\ (K \in \A'').
\end{eqnarray*}
Then for $D_0(\A):=\{\theta \in D(\A) \mid \theta(\alpha_{H_0})=0\}$, the restriction map 
$$
\pi:D_0(\A) \rightarrow D(\A'',m_0)
$$
is defined by taking a residue of $\alpha_{H_0}$. See \cite{Z} for details.  
\label{zr}
\end{define}

Next let us introduce the shift isomorphism, which will be generalized 
in Theorem \ref{univ} for $2$-arrangements. In this paragraph we assume that $\K=\R$. 
Let $\A$ be a Coxeter arrangement with the Coxeter 
group $W$. Put $R:=S^W=\R[P_1,\ldots,P_\ell]$ with 
homogeneous basic invariants $P_1,\ldots,P_\ell$ by Chevalley's theorem. 
Let $F$ be a quotient field of $S$. 
We may assume that 
$\deg P_\ell=h >\deg P_i\ (i\neq \ell)$, where $h$ is the Coxeter number of $W$. 
Let $D \in \Der(R):=\oplus_{i=1}^\ell R \cdot 
\partial_{P_i}$ be the invariant derivation of degree $-h+1$, called the \textbf{primitive 
derivation}. 
Then for 
the Euler derivation $\theta_E$ and the affine connection $\nabla$ defined by 
$$
\nabla_\theta (\sum_{i=1}^\ell f_i \partial_{x_i}):=
\sum_{i=1}^\ell \theta(f_i) \partial_{x_i}
$$
for $\theta \in \Der(F):= \Der (S) \otimes _S F$, the following shift isomorphism holds.

\begin{theorem}[\cite{AY2}, Theorem 2]
For $m:\A \rightarrow \{+1,0\}$ and $k \in \Z_{\ge 0}$, the $S$-morphism
$$
\Phi:D(\A,m) \rightarrow D(\A,2k+m)
$$
defined by 
$$
\Phi(\theta):=\nabla_\theta \nabla_D^{-k} \theta_E
$$
is an isomorphism.
\label{shift}
\end{theorem}

For the most generalized version of shift isomorphisms, see Theorem 0.7 in \cite{A4}. 
In the rest of this 
section assume that $\A$ is a $2$-arrangement in $\K^2$. 
Let us recall results in \cite{AN}.
For a multiarrangement $(\A,m)$ with $\exp (\A,m)=(d_1,d_2)$, recall that 
$$
\Delta(m)=|d_1 - d_2|.
$$
Then the \textbf{multiplicity lattice} $\Lambda$ and the subset $\Lambda'$ is defined by 
\begin{eqnarray*}
\Lambda:&=&\{m:\A \rightarrow \Z_{\ge 0}\},\\
\Lambda':&=&\{m \in \Lambda\mid \Delta(m) \neq 0\}.
\end{eqnarray*}
Then $\Delta$ is a function 
$$
\Delta:\Lambda \rightarrow \Z_{\ge 0}.
$$
Also, define 
\begin{eqnarray*}
\Lambda_K:&=&\{ m \in \Lambda \mid 
m(K) > \sum_{H \in \A \setminus \{K\}} m(H)\}\ (K \in \A),\\
\Lambda_0:&=&\Lambda' \setminus \cup_{H \in \A} \Lambda_H.
\end{eqnarray*}
Note that $\Lambda_0$ is denoted by $\Lambda'_{\phi}$ in \cite{AN}. Put 
$$
d(m,m'):=\sum_{H \in \A} | m(H)-m'(H)|\ (m,m' \in \Lambda).
$$
Then the following structure theorems hold for $\Lambda_0$.

\begin{lemma}[\cite{AN}, Lemma 4.2]
For $m_1,m_2 \in \Lambda$ such that 
$d(m_1,m_2)=2$, 
$m_1(H)=m_2(H)$ for $H \in \A \setminus \{H_0\}$ and 
$m_1(H_0)=m_2(H_0)+1$, it holds that 
$|\Delta(m_1)-\Delta(m_2)|=1$.
\label{one}
\end{lemma}

\begin{theorem}[\cite{AN}, Theorem 3.2]
Let $C \subset \Lambda_0$ be a maximal connected component of $\Lambda_0$. 
Then there exists the unique point $m \in C$, called the \textbf{peak point of $C$}, such that 
$$
\Delta(m) \ge \Delta(\mu)\ (\forall \mu \in C).
$$
Moreover, 
$$
C=\{\mu \in \Lambda' \mid d(m,\mu) < \Delta(m)\}
$$
and for $\mu \in C$, 
$$
\Delta(\mu)=\Delta(m)-d(m,\mu).
$$
\label{str}
\end{theorem}

The maximal connected component of $\Lambda_0$ in Theorem \ref{str} is 
just called a \textbf{(finite) component}, and 
$\Lambda_K$ an \textbf{infinite component}. 
Also, the following independency property holds. 

\begin{prop}
Let $m_1, m_2 \in \Lambda'$ such that $d(m_1,m_2)=2$, 
$\Delta(m_i)=1\ (i=1,2)$,  
$m_1(H)=m_2(H)$ for $H \in \A \setminus \{H_1,H_2\}$, 
$m_1(H_1)=m_2(H_1)+1$ and $m_1(H_2)+1=m_2(H_2)$. 
Assume that for two multiplicities 
\begin{eqnarray*}
\mu(H):&=&\max\{m_1(H),m_2(H)\}\ (H \in \A),\\
\mu'(H):&=&\mbox{mim}\{m_1(H),m_2(H)\}\ ( H \in \A),
\end{eqnarray*}
it holds that $\Delta(\mu)=\Delta(\mu')=0$. 
Then  the lower degree bases
$\theta_i$ of $D(\A,m_i)\ (i=1,2)$ 
are $S$-independent.
\label{next}
\end{prop}

\noindent
\textbf{Proof}. This is the special case of Lemma 4.17 in \cite{AN}. \owari

\begin{rem}
We always start a multiarrangement $(\A,m)$ such that $m:\A \rightarrow \Z_{>0}$. 
However, in the arguments in the rest of this article, it often happens that 
a new multiplicity $m'$ attains zero at some hyperplane $H \in \A$. However, as we have seen in the above, 
the theory in \cite{AN} is constructed in the multiplicity lattice 
$\Lambda=\{ m:\A \rightarrow \Z_{\ge 0}\}$. Hence there are no problems.
\end{rem}

\section{Proof of Theorem \ref{limit}}
In this section let us prove Theorem \ref{limit}.
\medskip

\noindent
\textbf{Proof of Theorem \ref{limit}}. 
Assume that $\Delta(m) >h-2$. 
We may assume that $\{x_1x_2=0\} \subset \A$, which may not 
be orthogonal. Then there exist derivations $\partial_{x_1}, \partial_{x_2}$ of 
degree zero such that $\langle \partial_{x_1}, \partial_{x_2} \rangle_{\K}=
\Der(S)$ and that $\partial_{x_i}(x_j)=\delta_{ij}$. 
Let $C$ be a connected component 
of the multiplicity lattice of $\A$ such that $m \in C$. Since 
$\Delta$ attains the maximum value at the peak point of $C$, we may assume that 
$m$ is the peak point of $C$. 
Let $\theta$ be the lower degree basis for $D(\A,m)$ with $\deg \theta=d$. Then 
$\exp(\A,m)=(d,d+\Delta(m))$. 
If we put $n_i:=m(x_i=0)\ (i=1,2),\ n_H:=m(H)\ (H \in \A)$ and   
\begin{eqnarray*}
\theta (x_j)&=& x_j^{n_j} f_j\ (j=1,2),\\
\theta (\alpha_H)&=& \alpha_H^{n_H} f_H\ (H \in \A \setminus \{x_1x_2=0\}),
\end{eqnarray*}
then, for $\{i,j\}=\{1,2\}$, 
\begin{eqnarray*}
\nabla_{\partial_{x_i}} \theta (x_j)&=& x_j^{n_j} \partial_{x_i}(f_j),\\
\nabla_{\partial_{x_i}} \theta (x_i)&=& x_i^{n_i-1} (n_i f_i+ x_i \partial_{x_i}(f_i)),\\
\nabla_{\partial_{x_i}} \theta (\alpha_H)&=& \alpha_H^{n_H-1}(n_H \partial_{x_i}(\alpha_H) 
f_H+\alpha_H
\partial_{x_i}(f_H))\ (H \in \A \setminus \{x_1x_2=0\}).
\end{eqnarray*}
Hence 
the derivation $\nabla_{\partial_{x_i}} \theta$ 
is belonging to $D(\A,m-m_i)$, where 
\[
m_i(H)=
\left\{
\begin{array}{rl}
0 & \mbox{if}\ H = \{x_j=0\},\\
1 & \mbox{if}\ H \neq \{x_j=0\}
\end{array}
\right.
\] 
with $\{i,j\}=\{1,2\}$.
Since $|m_i|=|\A|-1=h-1$ and $\Delta(m) > h-2$, it holds that 
$m-m_i \in C$ or on the boundary of $C$ for $i=1,2$. 
Since $\exp(\A,m-m_i)=(d,d+\Delta(m)-h+1)$ by Theorem \ref{str}, 
$d \le d+\Delta(m)-h+1$ and 
$\deg \nabla_{\partial_{x_i}} \theta=d-1\ (i=1,2)$, it holds that 
$$
\nabla_{\partial_{x_i}} \theta=0\ (i=1,2). 
$$
Since $\mbox{char}(\K)=0$, it holds that $\theta=a \partial_{x_1}+b \partial_{x_2}$ with 
$a, b \in \K$, which is a contradiction. 
\owari
\medskip

\begin{rem}
If $\mbox{char}(\K) =p >0$, 
then the statement in Theorem \ref{limit} does not hold. For example, 
assume that $p=2$ and consider a balanced multiarrangement $(\A,m)$ defined by 
$$
x_1^4x_2^4(x_1+x_2)^4=0.
$$
Then Theorem \ref{Saito} shows that 
$
x_1^4 \partial_{x_1}+x_2^4 \partial_{x_2}
$
and 
$x_1^8 \partial_{x_1}+x_2^8 \partial_{x_2}$ form a basis 
for $D(\A,m)$. Hence $\Delta(m)=4>|\A|-2=1$. 

On the other hand, the proof of Theorem \ref{limit} says that, if 
$\mbox{char}(\K)=p>0$, $m \in \Lambda_0$ is a peak point 
and $\Delta(m)>h-2$, then 
the degree of the lower degree basis for $D(\A,m)$ 
can be divided by $p$. See also \cite{Num}. 
\end{rem}





\section{Proof of Theorem \ref{univ}}
In this section we prove Theorem \ref{univ}.
\medskip

\noindent
\textbf{Proof of Theorem \ref{univ}}. 
The assumptions and Theorem \ref{limit} 
imply that $m_0$ is the peak point of 
some finite component $C \subset \Lambda_0$. 
We may assume that $\{x_1x_2=0\} \subset \A$ and take 
$\partial_{x_1}, \partial_{x_2}$ in the same way as those in the proof of Theorem 
\ref{limit}. 
By the argument in \cite{Y0}, 
$\Phi_0(\theta) \in D(\A,m_0+m-1)$ for $\theta \in 
D(\A,m)$. Put $d:=\deg \theta_0$. 
Then 
$$
|m_0+m-1|=(2d+h-2)+|m|-h=2d+|m|-2.
$$
On the other hand, for a basis $\{\theta_1,\theta_2\}$ for $D(\A,m)$, 
$$
\deg \nabla_{\theta_1} \theta_0+\deg \nabla_{\theta_2}\theta_0=
(|m|-2)+2d.
$$
Noting that all multiplicities on $\A$ are free, 
by Theorem \ref{Saito} and arguments in \cite{Y0}, 
it suffices to show that 
$\nabla_{\partial_{x_1}} \theta_0$ and $\nabla_{\partial_{x_2}} \theta_0$ are 
$S$-independent. Define two multiplicities $m_1$ and 
$m_2$ by 
\[
m_i(H)=
\left\{
\begin{array}{rl}
0 & \mbox{if}\ H \neq \{x_j=0\},\\
1 & \mbox{if}\ H =\{x_j=0\},
\end{array}
\right.
\] 
where $\{i,j\}=\{1,2\}$. Then $\nabla_{\partial_{x_i}} \theta_0 \in 
D(\A,m_0+m_i-1)$ by the same arguments as in the proof of Theorem \ref{limit}. 
Since $d(m_0,m_0+m_i-1)=h-1$, 
it holds that $\exp(\A,m_0+m_i-1)=
(\deg \theta_0-1, \deg \theta_0)$ by Lemma \ref{one} and 
Theorem \ref{str}. Note that if $\nabla_{\partial_{x_i}} \theta_0 \neq 0$, then it  
is of degree $\deg \theta_0-1$ and is the lower degree basis 
for $D(\A,m_0+m_i-1)$. Let us check that $\nabla_{\partial_{x_i}} \theta_0\neq 0$. 
We may assume that $i=1$ and suppose that 
$\nabla_{\partial_{x_1}}\theta_0=0$. Then 
$\theta_0=ax_2^d \partial_{x_1}+bx_2^d \partial_{x_2}\ (a,b \in \K)$, which 
is only tangent to one of following three arrangements of hyperplanes:
\begin{itemize}
\item[(1)]
$\{(bx_1-ax_2)x_2=0\}$ if $a \neq 0,\ b\neq0$.
\item[(2)]
$\{x_1x_2=0\}$ if $a =0,\ b\neq0$.
\item[(3)]
$\{x_2=0\}$ if $b=0$.
\end{itemize}
Any case contradicts $h>2$. Hence $\nabla_{\partial_{x_i}} \theta_0\neq 0\ (i=1,2)$ and 
both derivations are the lower degree basis of degree $d-1$ for $D(\A,m_0+m_i-1)$.

Note that 
$\Delta(m_0+m_i-1)=1$ and $\Delta(m_0+m_1+m_2-1)=0$. 
Hence, by Proposition \ref{next}, it suffices to show that 
$\Delta(m_0-1)=0$. 

Assume that $\Delta(m_0-1) \neq 0$. Then Lemma \ref{one} shows that 
$\Delta(m_0-1)=2$. 
Then the result in \cite{AN}, or the addition theorem in \cite{ATW2} says that 
$$
\nabla_{\partial_{x_1}} \theta_0=x_2  \theta', \ 
\nabla_{\partial_{x_2}} \theta_0=x_1  \theta'
$$
(up to scalars) for the lower degree basis $\theta' $ for $D(\A,m_0-1)$. 
Hence it holds that 
\begin{equation}
x_1\nabla_{\partial_{x_1}} \theta_0
=x_1x_2 \theta'
=
x_2\nabla_{\partial_{x_2}} \theta_0
\label{uts}
\end{equation}
up to scalars. 
Let us put $\theta_0=f\partial_{x_1}+g \partial_{x_2}$ and 
$f=\sum_{i=0}^d a_i x_1^{i} x_2^{d-i}$. Then (\ref{uts}) implies that 
$$
ia_i=(d-i)a_i\ (i=0,\ldots,d)
$$
up to a unique scalar. 
Note that $\Delta(m \equiv 1)=h-2$ and 
$|m \equiv 1 |=h>0$. Since 
there is nothing to prove when $m_0\equiv 1$, 
we may assume that 
$|m_0| -|\A| \ge 2$. Thus $d \ge 2$. 
Also, for $i \neq j$, 
$$
\displaystyle \frac{i}{d-i}=\displaystyle \frac{j}{d- j} \iff 
d=0.
$$ 
Hence 
$f$ is a monomial. The same argument to $g$ shows that 
$\theta_0$ is of the form
$$
\theta_0=a x_1^i x_2^{d-i} \partial_{x_1}+b x_1^i x_2^{d-i} \partial_{x_2}
= x_1^i x_2^{d-i}(a  \partial_{x_1}+b \partial_{x_2})
$$
with $1 \le i \le d-1$. Apparently 
$\theta_0$ is tangent to $x_1=0$ and $x_2=0$ with multiplicity $i$ and 
$d-i$ respectively, 
can be tangent to $bx_1-ax_2=0$ and 
tangent to no other hyperplanes. When $|\A|=3$ and 
$m_0-1 \in \Lambda_0$, the statement is true by \cite{Waka} 
and Proposition \ref{next}. 
If $|\A| \ge 4$ then $\theta_0$ 
cannot be in $D(\A,m_0)$. Hence $\Delta(m_0-1)=0$ and 
Proposition \ref{next} completes the proof. 
\owari
\medskip

%
%
%

Theorem \ref{univ} says that, if $\Delta(m)=h-2$, then 
$D(\A) \simeq D(\A,m)$ as $S$-modules. Since $\theta_0$ in Theorem \ref{univ} for 
$D(\A)$ is the Euler derivation $\theta_E$, we introduce the following definition. 

\begin{define}
The derivation $\theta_0$ in Theorem \ref{univ} is called to be  
the \textbf{$(\A,m)$-Euler derivation}. 
\label{E}
\end{define}

Obviously the Euler derivation $\theta_E$ is $(\A,1)$-Euler for all $2$-arrangements. 
Let us see the other examples below. 

\begin{example}
Let $(\A,m)$ be an $A_2$-type multiarrangement. By \cite{Waka} 
we know that $\Delta(m)=1$ if $|m|$ is odd and $m$ is balanced. Since $|\A|=3$, Theorem \ref{univ} 
shows that every lower degree basis for $D(\A,m)$ such that 
$|m|\equiv 1\ $(mod $2)$ and $m,m-1 \in \Lambda_0$ 
is 
$(\A,m)$-Euler. 
%
\end{example}

Historically a lot of $(\A,2k+1)$-Euler derivations have 
been constructed by using the invariant theory for Coxeter arrangements $\A$, 
see \cite{Y0} and \cite{AY2}. 
In these papers to prove the independency is an important part. By using 
Theorem \ref{univ}, we can give an another proof when $\ell=2$. 

\begin{cor}
Let $\A$ be a Coxeter arrangement in $\R^2$ corresponding to the Coxeter group $W$ and put 
$R:=S^W=\R[P_1,P_2]$ the invariant ring with basic invariants. Define $D:=\partial_{P_2} \in 
\Der(R)$ the primitive derivation. 
Then the derivation 
$E_{k}:=\nabla_{D}^{-k} \theta_E$ is 
$(\A,2k+1)$-Euler for $k \in \Z$.
\label{G2}
\end{cor}

\noindent
\textbf{Proof}. 
Assume that 
$k \ge 0$. Then 
$\deg E_k=hk+1$ with $h=|\A|$. 
Also, $|(\A,2k+1)|=h(2k+1)$ and a constant multiplicity is balanced. 
Since 
$E_k \in D(\A,2k+1)$ by \cite{AT}, Theorem \ref{univ} completes 
the proof. 

Assume that $k<0$. Note that the same theory in \cite{AN} holds true for 
$-\Lambda:=\{m:\A \rightarrow \Z_{\le 0}\}$. Then, noting that 
$\nabla_\theta(\omega) \wedge d\alpha =
\nabla_\theta (\omega \wedge d \alpha)$ for 
a rational differential form $\omega$ and 
$\alpha \in V^*$, the same argument as the above 
completes the proof. 
\owari
\medskip

The proof of Theorem \ref{univ} implies the following. 

\begin{prop}
Let $(\A,m_0)$ be a $2$-multiarrangement with $|\A|=h \ge 4$, $m_0 \in \Lambda_0$ and 
$\Delta(m_0)=h-2$. Then $m_0-1$ is balanced. 
\end{prop}

\noindent
\textbf{Proof}. Assume not. Then 
it is obvious that $\nabla_{\partial_{x_1}} \theta_0$ 
and $\nabla_{\partial_{x_2}} \theta_0$ in the proof of Theorem \ref{univ} 
are $S$-dependent since 
they are in the same infinite component and that $\Delta(m_0-1)=2$. 
\owari

\section{Freeness condition for $3$-arrangements}

Before the proof of Theorem \ref{fc} let us recall 
one of Yoshinaga's freeness criterions. 

\begin{theorem}[\cite{Y3}, Theorem 3.2]
Let $\A$ be a central $3$-arrangement with the infinite hyperplane $H_0 \in \A$. 
Let $(\A \cap H_0,m_0)$ be the Ziegler restriction onto $H_0$ with 
$\exp(\A \cap H_0, m_0)=(d_1,d_2)$. Also put 
$$
\chi(\A,t)=(t-1)(t^2-c_1t+c_2).
$$
Then, for the Ziegler restriction map 
$$
\pi:D_0(\A) \rightarrow D(\A \cap H_0,m_0),
$$
it holds that 
$$
\dim \coker \pi =
c_2-d_1d_2 \ge 0.
$$
In particular, $\A$ is free with $\exp(\A)=(1,d_1,d_2)$ if and only if 
$c_2=d_1d_2$.
\label{Y}
\end{theorem}

The following is an immediate consequence of Theorem \ref{Y} and well-known, see 
\cite{WY} for example. Here 
we give a proof. 

\begin{prop}
The Terao conjecture is true for a $3$-arrangement 
$\A$ such that its Ziegler restriction is not balanced. 
\label{nb}
\end{prop}

\noindent
\textbf{Proof}. 
Fix an infinite 
hyperplane $H_0 \in \A$ and put $(\A'',m)$ the Ziegler restriction of $\A$ onto $H_0$ 
such that $m(K) > \sum_{H \in \A'' \setminus \{K\}}m(H)$ for some $K \in \A''$. 
We may assume that 
$\alpha_K=x_1$. Then 
$$
\prod_{H \in \A'' \setminus \{K\}} \alpha_H^{m(H)} \partial_{x_2}
$$
is the lower degree basis for $D(\A'',m)$. Hence 
$\exp(\A,m)=(m(K),|m|-m(K))$ and Theorem \ref{Y} says that 
$\A$ is free if $\chi(\A,t)=(t-1)(t-m(K))(t-|m|+m(K))$, which is combinatorial. \owari
\medskip

Now let us consider the balanced cases. 
Let us recall a notation. 
Let $\overline{\A}$ be an affine $2$-arrangement, $c\overline{\A}$ its coning 
with the 
infinite hyperplane $H_0 \in c \overline{\A}$. Put $|\overline{\A}|=k$ and 
$\chi(\overline{\A},t)=(t^2-kt+c_2)$. Also, let $(\A,m)$ be 
the Ziegler restriction of $c \overline{\A}$ onto $H_0$ with $|\A|=h>2$. 
We say that 
an affine simple $2$-arrangement $\overline{\A}$ is \textbf{balanced} if 
$m(K) \le \sum_{H \in \A \setminus \{K\}} m(H)$ for any $K \in \A$. Now 
we have prepared for the proof of Theorem \ref{fc}. 
\medskip


\noindent
\textbf{Proof of Theorem \ref{fc}}. 
First assume that $c_2=d\times (d+h-2)$ and 
$k=2d+h-2$. Put 
$\exp(\A,m)=(d_1,d_2)$. Note that $d_1+d_2=2d+h-2=k$. 
Since $m$ is balanced, 
Theorem \ref{limit} implies that 
$$
d_1d_2 \ge d(d+h-2).
$$
Also, Theorem \ref{Y} implies that 
$$
c_2=d(d+h-2) \ge d_1d_2.
$$
If $\exp(\A,m) \neq (d,d+h-2)$, then 
$d < \min\{d_1,d_2\}$ and 
$\max\{d_1,d_2\} < d+h-2$ by Theorem \ref{limit}. Hence 
$$
c_2=d(d+h-2) \ge  d_1d_2 > d(d+h-2), 
$$
which is a contradiction. Hence 
$\exp(\A,m) = (d,d+h-2)$, and Theorem \ref{Y} 
completes the proof. 

Second assume that 
$c_2=d\times (d+h-3)$ and 
$k=2d+h-3=d_1+d_2$. In this case, 
$k-(h-2)$ is an odd number. Hence 
$\Delta(m)=|d_1-d_2| \le h-3$. By the same arguments as the above, 
$$
c_2=d(d+h-3) \ge d_1d_2 \ge d(d+h-3),
$$
which, combined with $h=|\A|>2$, completes the proof. \owari
\medskip

\noindent
\textbf{Proof of Corollary \ref{TC3}}. 
Immediate by Theorem \ref{fc}. 
\owari
\medskip

Theorem \ref{limit} has a lot of applications on 
the characteristic polynomials, freeness and chambers of $3$-arrangements as 
follows. 

\begin{theorem}
In the above notation, 
assume that $\overline{\A}$ is balanced and 
$\chi(\overline{\A},t)=(t-a)(t-b)$ with $a \le b$. 
\begin{itemize}
\item[(1)]
If $k=a+b=2d+h-2$ for some integer $d$, 
then $d \le a \le b \le d+h-2$. 
\item[(2)]
If $k=a+b=2d+h-3$ for some integer $d$, 
then $d \le a \le b \le d+h-3$.
\end{itemize}
\label{rest}
\end{theorem}

\noindent
\textbf{Proof}.  
Since the proof is the same, we only prove (1). 
By Theorem \ref{Y} it holds that 
$$
ab \ge d_1d_2
$$
with $\exp(\A,m)=(d_1,d_2)$. Note that $d_1+d_2=a+b$. By Theorem \ref{limit}, 
$|d_1-d_2| \le h-2$. Hence, if 
$b-a >h-2$, then $ab-d_1 d_2 <0$, which is a contradiction. \owari
\medskip

\begin{theorem}
In the above notation, 
assume that $\overline{\A}$ is balanced and $\K=\R$. Let $\mbox{ch}(\overline{\A})$ 
be the set of connected components of $\R^2 \setminus \cup_{H \in \overline{\A}}  H$. 
\begin{itemize}
\item[(1)]
If $k=2d+h-2$ for some integer $d$, 
then $c_2 \ge d(d+h-2)$ and 
$$
|\mbox{ch}(\overline{\A})| \ge 1+k+d(d+h-2).
$$
In particular, the equation holds only if 
$c\overline{\A}$ is free. 

\item[(2)]
If $k=2d+h-3$ for some integer $d$, 
then $c_2 \ge d(d+h-3)$ and 
$$
|\mbox{ch}(\overline{\A})| \ge 1+k+d(d+h-3).
$$
In particular, the equation holds only if 
$c\overline{\A}$ is free. 

\end{itemize}
\label{rest2}
\end{theorem}

\noindent
\textbf{Proof}.
The same as that of Theorem \ref{rest}. \owari
\medskip

These results say that, if $\overline{\A}$ is balanced, then 
the characteristic polynomial $\chi(\overline{\A},t)$ 
is irreducible, or splits with a restricted splitting type 
seen in the above. Also, the freeness of these arrangements are determined by 
the intersection lattice, or more explicitly, by the characteristic polynomial. 

The following can be proved by using results in \cite{Waka}. We give an another proof here. 

\begin{cor}
Let $c \overline{\A}$ be a $3$-arrangement with the infinite hyperplane $H_0 \in c \overline{\A}$ and 
$(\A,m)$ the Ziegler restriction of $c \overline{\A}$ onto $H_0$. 
Assume that $|\A|=3$. Then the freeness of $c \overline{\A}$ depends only on 
$L(\overline{\A})$. In particular, so is that of the deformation of 
the Coxeter arrangement of type $A_2$. 
\label{A2}
\end{cor}

\noindent
\textbf{Proof}. 
If $(\A,m)$ is not balanced, then Proposition \ref{nb} completes the proof. 
Assume that $(\A,m)$ is a balanced Coxeter multiarrangement of type $A_2$. 
Since $|\A|-2=1$, Theorem \ref{fc} completes the proof. \owari
\medskip

\begin{example}
Let $\A$ be the Coxeter arrangement of type $B_2$ and $\overline{\A}$ its 
deformation as in \cite{A2}. The freeness of such deformations have not yet 
classified, nor have the exponents of the multiarrangements $(\A,m)$. 
Some of the freeness of $c \overline{\A}$ was classified in Propositions 
2.3 and 2.4 in \cite{A2}. 
They were proved by using the addition theorem. If we use Theorem \ref{fc}, 
the explicit formula of the Poincar\'e polynomial (\cite{A2}, Lemma 2.1) is enough to 
show the freeness. In other words, Theorem \ref{rest} says that if $\overline{\A}$ is 
balanced and $c\overline{\A}$ splits, then 
it is of the form 
\begin{eqnarray*}
\chi(\overline{\A},t)&=&(t-d)^2,\\
\chi(\overline{\A},t)&=&(t-d)(t-d+1)\ \mbox{or}\\
\chi(\overline{\A},t)&=&(t-d)(t-d+2).
\end{eqnarray*}
Then Theorem \ref{fc} says that $c \overline{\A}$ is free if $\chi(\overline{\A},t)$ splits 
into the form of the second and third types in the above. 
\end{example}

More generally, when $|\A|=4$ and not necessarily of type $B_2$, the following holds. 

\begin{cor}
Let $\A$ be a central $3$-arrangement such that 
$|\A \cap H_0|=4$ for some $H_0 \in \A$. If 
$|\A|$ is even, then the freeness of $\A$ depends only on 
$L(\A)$.
\label{four}
\end{cor}

\noindent
\textbf{Proof}. 
If the deconing $d\A:=(\A \setminus \{H_0\})|_{\alpha_{H_0}=1}$ (a 
converse operation of the coning) 
is not balanced, then Proposition \ref{nb} completes the proof. 
Assume that $d\A$ is balanced. 
In this assumption, the splitting type of the characteris	tic polynomial is 
always 
$$
\chi(\A,t)=(t-1)(t-d)(t-d-1)
$$
by Theorem \ref{rest}. Hence Theorem \ref{fc} completes the proof. \owari
\medskip

 \vspace{5mm}

\end{document}